\documentclass{amsart}
\usepackage[dvips]{graphicx}

\newtheorem{theorem}{Theorem}[section]
\newtheorem{lemma}[theorem]{Lemma}

\newtheorem{corollary}[theorem]{Corollary}

\theoremstyle{definition}
\newtheorem{definition}[theorem]{Definition}
\newtheorem{remark}[theorem]{Remark}

\numberwithin{equation}{section}

\newcommand{\zz}{\mathbb{Z}}
\newcommand{\cc}{\mathbb{C}}

\begin{document}

\title[Homologous Non-isotopic Symplectic Tori
in a $K3$--Surface]{Homologous Non-isotopic Symplectic Tori\\
in a \mbox{\boldmath $K3
\frac{\hspace{4pt}}{\hspace{4pt}}$}Surface}

\author{Tolga Etg\"u}
\address{Department of Mathematics and Statistics, McMaster University,
Hamilton, Ontario L8S 4K1, Canada} \email{etgut@math.mcmaster.ca}

\author{B. Doug Park}
\address{Department of Pure Mathematics, University of Waterloo, Waterloo,
Ontario, N2L 3G1, Canada} \email{bdpark@math.uwaterloo.ca}

\thanks{B.D. Park was partially supported by NSERC and CFI grants.}

\date{April 4, 2003.  Revised on July 14, 2003}

\begin{abstract}
For each member of an infinite family of homology classes in the
$K3$--surface $E(2)$, we construct infinitely many non-isotopic
symplectic tori representing this homology class. This family has
an infinite subset of primitive classes. We also explain how these
tori can be non-isotopically embedded as homologous symplectic
submanifolds in many other symplectic 4-manifolds including the
elliptic surfaces $E(n)$ for $n>2$.
\end{abstract}

\maketitle


\section{Introduction}

A homology class in a complex surface is represented by at most
finitely many complex curves up to smooth isotopy. In contrast,
there are examples of symplectic 4-manifolds admitting infinite
families of homologous but non-isotopic symplectic submanifolds
(see e.g. \cite{ep1}, \cite{fs:non-isotopic},
\cite{vidussi:non-isotopic}).  For example, in \cite{ep1}, we
constructed infinitely many homologous, non-isotopic symplectic
tori representing the divisible homology class $q[F]$, for each $q
\ge 2$, where $F$\/ is a regular fiber of a simply-connected
elliptic surface $E(n)$\/ with no multiple fibers. In this paper
we construct such infinite families in the homology class
$q[F]+m[R]$, for any pair of positive integers $(q,m)\neq(1,1)$,
where $[R]$ is the homology class of a rim torus in $E(n)$ with
$n\geq 2$.  In particular, we get non-isotopic tori in infinitely
many \emph{primitive}\/ homology classes. Unfortunately, primitive
classes in $E(1)$ seem to be still out of our reach at the moment.
Examples of tori representing primitive homology classes in
symplectic 4-manifolds homeomorphic to $E(1)$ are given in
\cite{ep:E(1)_K} and \cite{vidussi:E(1)_K}.

A significant difference between the construction we give here and
the examples in \cite{ep1}, \cite{fs:non-isotopic} and
\cite{vidussi:non-isotopic}  is that the tori here are not
obtained by braiding of parallel copies of the same symplectic
surface (a regular fiber $F$\/ of an elliptic fibration) in the
sense of \cite{adk}, but rather using parallel copies of two
different symplectic surfaces ($F$ and a rim torus $R$).  In fact,
$R$\/ is Lagrangian with respect to the symplectic form on $E(n)$
induced by the elliptic fibration.  In some cases we need to use a
small perturbation of this symplectic form with respect to which
$R$\/ becomes symplectic.

As a consequence of our calculations,
we are able to distinguish the tori we construct not only
up to smooth isotopy
but also up to self-diffeomorphisms of the ambient 4-manifold.
We should also note that, just like our earlier result in \cite{ep1},
the construction here extends
to a more general class of symplectic 4-manifolds
(see Theorem~\ref{thm:generalization}).
In the sequel \cite{ep2}, we construct families of homologous
non-isotopic Lagrangian tori using different methods.

In the next section, we state our main result,
Theorem~\ref{theorem:main}, after a brief review of some basic
facts about the complex elliptic surface $E(2)$, which is a
$K3$--surface. (For more details on the topology of $E(2)$ and
other elliptic surfaces, we refer to the excellent book
\cite{gs}.) In Sections \ref{sec:link surgery}--\ref{sec:tori in
E(2)}, we explain two general constructions which utilize braids
to give symplectic tori in $E(2)$ within a prescribed homology
class. In Section \ref{sec:alexander}, we apply these
constructions using particular set of braids which are suitable
for certain Seiberg-Witten invariant calculations. In Section
\ref{sec:sw}, we explain how these invariants distinguish the
symplectic tori up to isotopy. In the last section, we discuss
some possible generalizations of Theorem~\ref{theorem:main} to
other symplectic 4-manifolds.


\section{Topology of the $K3$--Surface $E(2)$ and the Main Result}

$E(2)$ is simply-connected. The intersection form of $E(2)$ is $2
E_8 \oplus 3H$, where $E_8$ is a unimodular negative definite
8$\times$8 matrix and $H:= \bigl(
\begin{smallmatrix}
  0 & 1 \\
  1 & 0
\end{smallmatrix} \bigr)$.
Let $[F],[S]$ denote the homology classes of a regular fiber and a
section of an elliptic fibration $f\! :\! E(2)\rightarrow
\mathbb{CP}^1$, respectively. They correspond to one summand of
$H$ in the intersection form. $E(2)$ is the fiber sum,
\[ E(2) = E(1)\#_F E(1) =
[E(1)\setminus \nu F] \cup_{\varphi} [E(1)\setminus \nu F],
\]
where a tubular neighborhood $\nu F$\/ is canonically identified
with the Cartesian product $F\times D^2$, and the gluing
diffeomorphism $\varphi :\partial (\nu F) \rightarrow
\partial (\nu F)$\/ identifies the fibers and is the complex
conjugation on the boundary of any normal disk, $\{{\rm
point}\}\times D^2$.  We fix a Cartesian product decomposition\/
$F= C_1 \times C_2$, where each
$C_j \cong S^1$.  Let $R_1 = C_1 \times
\partial D^2$, $R_2 = C_2\times \partial D^2 \subset E(2)$.
$R_i$ are called \emph{rim tori}.
Each circle $C_i$ bounds a disk in both copies of $[E(1)\setminus
\nu F]$ and gluing together the disks from both sides, we get a
sphere of self-intersection $-2$\/ in $E(2)$, which we denote by
$D_i$.  The remaining two $H\cong \bigl(
\begin{smallmatrix}
  0 & 1 \\
  1 & -2
\end{smallmatrix} \bigr) $\/ summands are generated by the homology bases
$\{[R_1], [D_2]\}$ and $\{[R_2],[D_1]\}$.
Our first result is the following.

\begin{theorem}\label{theorem:main}
For any pair of positive integers\/ $(q,m) \neq (1,1)$
there exists an infinite family of pairwise
non-isotopic symplectic tori representing the homology class\/
$q[F] + m[R_i]$\/ $(i=1$ or\/ $2)$ of an elliptic surface $E(2)$,
where\/ $[F]$\/ is the homology class of the fiber, and $[R_1]$
and $[R_2]$ are the homology classes of the rim tori.
\end{theorem}

\begin{remark}
Note that when $q$ and $m$ are relatively prime we obtain an
infinite family of pairwise non-isotopic symplectic tori
representing the same \emph{primitive}\/ homology class in $E(2)$.
\end{remark}

The proof of Theorem~\ref{theorem:main} is spread out over the
next few sections.


\section{Link Surgery}
\label{sec:link surgery}

We review the generalization of the link surgery construction of
Fintushel and Stern \cite{fs:knots} by Vidussi
\cite{vidussi:smooth}.  For an $n$-component link $L\subset S^3$,
choose an ordered homology basis of oriented simple curves
$\{\alpha_i, \beta_i \}_{i=1}^{n}$ such that $\alpha_i$\/ and\/
$\beta_i$ lie in the $i$-th boundary component of the link
exterior and the intersection number of $\alpha_i$ and $\beta_i$
is 1. Let $X_i$ ($i=1,\dots, n$) be a 4-manifold containing a
2-dimensional torus submanifold $F_i$ of self-intersection $0$.
Choose a Cartesian product decomposition $F_i = C_1^{i} \times
C_2^{i}$, where each $C^i_j \cong S^1$ ($j=1,2$)\/ is an embedded
circle in $X_i$.

\begin{definition}\label{def:data}
The ordered collection\/ $\mathfrak{D} := (\{(\alpha_i, \beta_i)
\}_{i=1}^{n}, \{X_i, F_i=  C_1^{i} \times C_2^{i} \}_{i=1}^{n})$\/
is called a \emph{link surgery gluing data}\/ for an $n$-component
link $L$. We define the \emph{link surgery manifold corresponding
to} $\mathfrak{D}$ to be the closed $4$-manifold
\[
L(\mathfrak{D}) \: :=\; [\coprod_{i=1}^{n} X_i\setminus\nu
F_i]\hspace{-20pt}\bigcup_{F_i\times\partial D^2=(S^1\times
\alpha_i)\times\beta_i}\hspace{-20pt} [S^1\times(S^3\setminus \nu
L)]\, ,
\]
where $\nu$\/ denotes the tubular neighborhoods.  Here, the gluing
diffeomorphisms between the boundary 3-tori identify the torus
$F_i =  C_1^{i} \times C_2^{i}$\/ of $X_i$\/ with\/ $S^1\times
\alpha_i$\/ factor-wise, and act as the complex conjugation on the
last remaining\/ $S^1$ factor.
\end{definition}

\begin{figure}[!ht]
\begin{center}
\includegraphics[scale=.5]{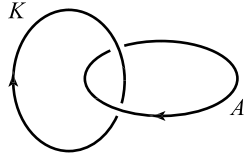}
\end{center}
\caption{Hopf link $L$} \label{fig:hopf}
\end{figure}

\begin{lemma}\label{lemma:E(2)}
Let\/ $L\subset S^3$ be the Hopf link in Figure\/
$\ref{fig:hopf}$.
For the link surgery gluing data\/
\begin{equation}\label{eq:data}
\mathfrak{D}=\big( \{ (\mu(K),\lambda(K)), (\lambda(A),-\mu(A) \},
\{ E(1), F=C_1\times C_2 \}_{i=1}^2 \big),
\end{equation}
we obtain\/ $L(\mathfrak{D}) \cong E(2)$.  Here, $\mu(K)$ and $\lambda(K)$ denote
the meridian and the longitude
of the knot $K$, respectively.
\end{lemma}

\begin{proof}
Note that the exterior of the Hopf link\/ $(S^3 \setminus \nu
L)$\/ is diffeomorphic to $S^1 \times \mathbb{A}$, where
$\mathbb{A}\cong S^1\times [0,1]$\/ is an annulus.  Hence there is
a diffeomorphism between the cylinder $\partial(\nu F)\times[0,1]
= T^2 \times
\partial D^2 \times [0,1]$\/ and the Cartesian product\/ $[S^1 \times (S^3
\setminus \nu L)]$.  We can easily check that our link surgery
gluing data is consistent with the fiber sum construction, and
gives
\[
L(\mathfrak{D}) = [E(1) \setminus \nu F] \cup [\partial(\nu
F)\times[0,1]\hspace{1pt} ] \cup [E(1)\setminus \nu F] \cong E(2)
.
\qedhere
\]
\end{proof}


\section{Two Symplectic Forms on the Cylinder\/ $\partial(\nu F)\times[0,1]$}
\label{sec:cylinder}

Let $M:=(S^3\setminus \nu L)$\/ denote the complement of the
tubular neighborhood of a 2-component Hopf link $L$\/ in $S^3$. We
saw that $M$ is diffeomorphic to a solid torus minus a thickened
core, i.e. $M \cong S^1 \times \mathbb{A}(r_0,r_1)$, where
$\mathbb{A}(r_0,r_1)=\{z\in \cc \: :\: r_0\leq |z| \leq r_1 \}$.
(In Figure~\ref{fig:annulus} the core is represented by the
darkened circle wherein you have no ``pineapple''.)  Let
$(r,\theta)$ be the polar coordinates on the annulus
$\mathbb{A}(r_0,r_1)$ with $-\pi < \theta \leq \pi$. Let
$(y,r,\theta)$ be the coordinate system on $M= S^1 \times
\mathbb{A}(r_0,r_1)$, where $y$ denotes the angular coordinate on
the $S^1$ factor ($-\pi< y \leq \pi$). For the sake of
concreteness, let us assume from now on that $r_1 = r_0 +1$.

Now define a 4-manifold with boundary $Y:=S^1 \times M \cong
[\partial(\nu F)\times[0,1] \hspace{1pt} ]$, and let $x$ be the
angular coordinate on the first $S^1$ factor ($-\pi < x \leq
\pi$).  To distinguish this $S^1$ factor with coordinate $x$ from
the $S^1$ factor in $M$\/ with coordinate $y$, we will denote them
by $S^1_x$ and $S^1_y$, respectively.

\subsection{First Family of Tori}
Our first symplectic form on $Y$\/ will be
\begin{equation}\label{eq:symplectic form}
dx \wedge dy \,+\, r dr\wedge d\theta\, = \,
\left.\omega_{f}\right|_{\partial(\nu F)\times[0,1] } ,
\end{equation}
where $\omega_{f}$ is the symplectic form on $E(2)$ coming from
the elliptic fibration $f$\/ (see Section~\ref{sec:tori in E(2)}).
Now let $B$\/ be a $q$-strand braid whose closure $\hat{B}$\/ is a
single-component link, i.e. a knot.  It is not hard to embed
$\hat{B}$ into the link exterior $M$\/ such that $S^1_x \times
\hat{B} \subset Y$\/ is a symplectic submanifold with respect to
the symplectic form (\ref{eq:symplectic form}).  We choose a
particular family of embeddings shown in Figure~\ref{fig:altenate
embedding}. Here we require the linking numbers to be
$lk(\hat{B},K)=q\,$ and $lk(\hat{B},A)=m$.

\begin{figure}[!ht]
\begin{center}
\includegraphics[scale=.9]{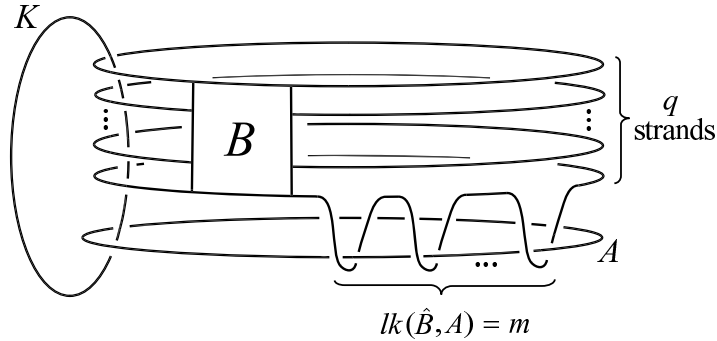}
\end{center}
\caption{A family of embeddings of $\hat{B}$\/ into\/
$(S^3\setminus\nu L)$} \label{fig:altenate embedding}
\end{figure}

Let us denote this family of embeddings by $\phi_{q,m} :\hat{B}
\rightarrow M$.  Then we have the following.

\begin{lemma}
For every pair of integers\/ $q\geq 2$\/ and\/ $m\geq 1$, the
torus\/ $S^1_x\times \phi_{q,m}(\hat{B})$\/ is a symplectic
submanifold of\/ $Y$ with respect to the symplectic form\/
$(\ref{eq:symplectic form})$.
\end{lemma}

\begin{proof}
We easily see that the restriction of the symplectic form
(\ref{eq:symplectic form}) to\/ $S^1_x\times\phi_{q,m}(\hat{B})$\/
is going to be just the restriction of\/ $dx\wedge dy$, which does
not vanish if we can arrange to have $dy\neq 0$\/ on the curve
$\phi_{q,m}(\hat{B})$. But this is always possible since we can
embed $\hat{B}$\/ in such a way that it is transverse to every
annulus of the form, $\{ {\rm point} \}\times
\mathbb{A}(r_0,r_1)$, inside $M$.
\end{proof}

\subsection{Second Family of Tori}
Our second symplectic form on $Y$\/ will be
\begin{equation}\label{def:omega}
\omega_s := dx \wedge (dy  + s \cdot d \theta) + r dr\wedge
d\theta\, ,
\end{equation}
where\/ $s>0$\/ is a sufficiently small real constant to be
determined later (see Section~\ref{sec:tori in E(2)}).  We easily
check that\/ $d\omega_s = 0$, and
\begin{equation*}
\omega_s\wedge\omega_s \,=\, 2 r \,dx\wedge dy\wedge dr\wedge
d\theta \,\neq\, 0\, .
\end{equation*}

\begin{figure}[!ht]
\begin{center}
\includegraphics[scale=.7]{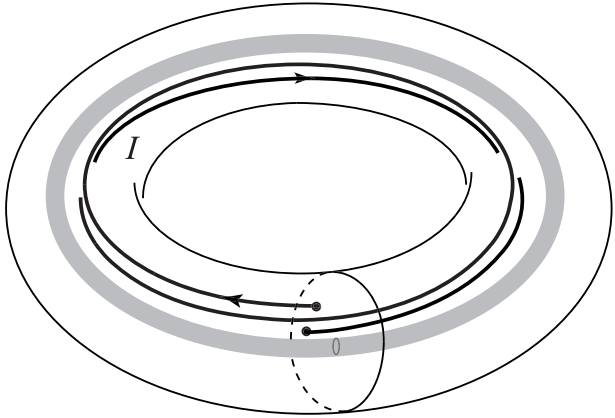}
\end{center}
\caption{An embedding of $\hat{B}$\/ into $M$\/ with\/ $m=2$}
\label{fig:torus}
\end{figure}

\begin{figure}[!ht]
\begin{center}
\includegraphics[scale=.7]{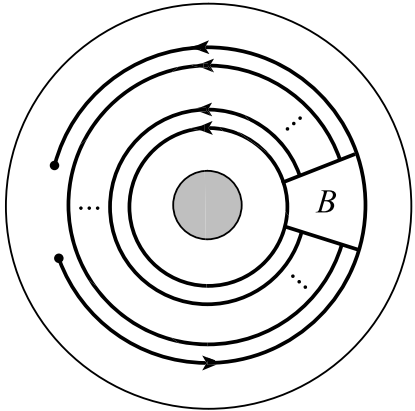}
\end{center}
\caption{Bird's eye view of the ``pineapple slice''}
\label{fig:annulus}
\end{figure}

Let $B$\/ be a $q$-strand braid as before.  We describe an
alternative way to embed the closure $\hat{B}$ into $M$.  (See
Figures~\ref{fig:torus} and \ref{fig:annulus}.)  Except for a
single connected arc $I$, the closed braid $\hat{B}$ lies inside a
thin ``pineapple slice'' of height $2 \varepsilon$, $\{(y
,r,\theta)\: : \: -\varepsilon\leq y \leq \varepsilon \}$.  The
remaining single arc $I$\/ traverses $m$ times around the solid
torus (minus core) in the positive $y$-direction ($m\geq 1$). Away
from the crossings and $I$, we require the circular arcs of
$\hat{B}$ to lie on a fixed level annulus, $\mathbb{A}_0
:=\{\,(y,r,\theta) \in M \: :\: y=0, \; r_0+\frac{1}{4} \leq r\leq
r_0+\frac{3}{4} \,\}$.  Note that the linking numbers are now
$\,lk(\hat{B},K)=m\,$ and $\,lk(\hat{B},A)=q$.  Essentially what
we are doing differently this time is reversing the roles of $K$\/
and $A$\/ in our first family of embeddings $\phi_{q,m}$ above.
(Compare Figures \ref{fig:altenate embedding} and
\ref{fig:braidform}.)

Obviously, $S^1_x \times (S^1_y \times \{{\rm point}\})$\/ is a
symplectic torus in $Y$ with respect to $\omega_s$.  We show that
$\hat{B}$ can be embedded into $M$\/ so that\/ $S^1_x \times
\hat{B}$ is also a symplectic torus in $(Y,\omega_s)$.  The
crucial condition is that the restriction of the 1-form\/ $\eta_s
:= dy + s \cdot d\theta \in \Omega^1(M)$\/ has a fixed sign over
the curve $\hat{B}$.

First orient the curve $\hat{B}$ as in Figures \ref{fig:torus} and
\ref{fig:annulus}.  Let\/ $\gamma: [0,\ell] \rightarrow \hat{B}$
be a parametrization of $\hat{B}$\/ by arc-length.  On the arc
$I$, we may arrange to have
\[
\langle dy , \dot{\gamma} \rangle \approx 1\, \quad {\rm and}
\quad \langle d\theta, \dot{\gamma}\rangle \geq 0
\]
as we traverse along the arc $I$\/ in the direction of the chosen
orientation.  This is possible because we can always embed $I$\/
so that $I$\/ is very close to being parallel to the (removed)
core of the solid torus.  Hence $\langle
\eta_s,\dot{\gamma}\rangle>0$, i.e. the 1-form $\eta_s$ is always
positive on $I$ in the chosen direction.

Next note that $\langle d\theta,\dot{\gamma}\rangle =1$, and
$\,\langle dr,\dot{\gamma}\rangle =\langle
dy,\dot{\gamma}\rangle=0$, away from the crossings in
$\mathbb{A}_0$.  Hence the restriction of $\eta_s$ is positive on
$\hat{B}\cap \mathbb{A}_0$, away from the crossings.

At a crossing in $B$, both $r$ and $y$ vary, so we need to draw
the braid such that
\begin{equation}\label{ineq:slopes}
\left|\frac{dy}{d\theta}\right| \,=\, \left|
\frac{dy/dt}{d\theta/dt} \right| \,<\,  s \: .
\end{equation}
Since we always have\/ $\langle d\theta, \dot{\gamma} \rangle =
d(\theta\circ\gamma)/dt = d\theta/dt
>0$\/ at any crossing,
an easy triangle inequality argument shows that\/ $\langle
\eta_s,\dot{\gamma}\rangle>0$\/ at every crossing.

In other words, we need to embed the braid $B$\/ so that every
pair of crossing arcs looks very short in terms of height $y$.
More precisely, we need to ensure that, as we traverse along the
crossing arcs in counter-clockwise direction, the angle $\theta$\/
is changing at a much faster rate than the rate of change for the
height $y$.

\begin{figure}[!ht]
\begin{center}
\includegraphics[scale=.7]{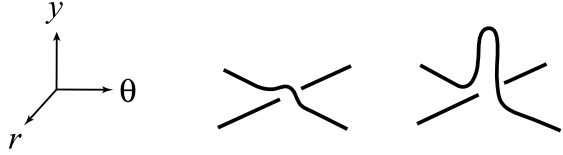}
\end{center}
\caption{A ``good'' crossing and a ``bad'' crossing.}
\label{fig:crossing}
\end{figure}

In Figure~\ref{fig:crossing}, the left crossing is short-looking
and hence ``good'', while the right crossing is something that we
must avoid.  To satisfy (\ref{ineq:slopes}) for small values of
$s$, we will have to embed the crossing arcs of $B$\/ very flat.
However note that there is no limitation on the number of
crossings or the number of strands allowed.

\begin{figure}[!ht]
\begin{center}
\includegraphics[scale=.7]{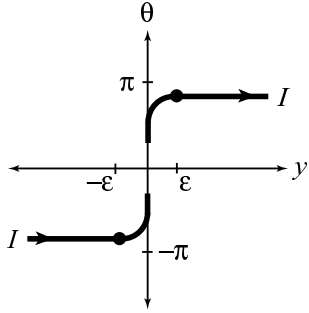}
\end{center}
\caption{Smooth corners at\/ $r={\rm constant}$} \label{fig:corner}
\end{figure}

Finally we need to verify that $\eta_s$\/ is positive on the two
``corners'' (which are represented by the two black dots in
Figures \ref{fig:torus} and \ref{fig:annulus}) where the arc $I$\/
is being attached to the rest of $\hat{B}$.  Note that we can
always assume that $r$ is constant on these two attaching portions
of $\hat{B}$. We can easily smooth out the corners such that
$\langle d\theta , \dot{\gamma}\rangle \geq 0$, $\langle dy ,
\dot{\gamma}\rangle \geq 0$, and the two quantities do not
simultaneously vanish (see Figure~\ref{fig:corner}).  Hence the
restriction of\/ $dy + s \cdot d\theta$\/ to the two corners is
strictly positive on the velocity vector $\dot{\gamma}$.

We conclude that $\eta_s$\/ restricts to some positive function
multiple of the orientation 1-form on $\hat{B}$.  Hence\/
$\omega_s|_p = dx\wedge \eta_s |_p \neq 0$, for every point\/
$p=(x,\gamma(t)) \in S^1_x \times \hat{B}$. Let us denote this
family of embeddings we constructed by\/ $\psi_{q,m} :\hat{B}
\rightarrow M$. In summary, we have the following.

\begin{lemma}\label{lemma:symplectic}
$\omega_s$\/ is a symplectic form on\/ $Y\!=\partial(\nu
F)\times[0,1]$ with respect to which the torus\/ $S^1_x \times
\psi_{q,m}(\hat{B})$ is a symplectic submanifold for every pair of
integers\/ $q \geq 2$\/ and\/ $m \geq 1$.
\end{lemma}


\section{Two Families of Homologous Symplectic Tori in $E(2)$}
\label{sec:tori in E(2)}

\begin{lemma}\label{lemma:symplectic form on E(2)}
There exists a symplectic\/ $2$-form\/ $\omega_f$ on $E(2)$, with
respect to which  the surfaces\/ $F$\/ and\/ $S$\/ are symplectic
and\/ $R_1$\/ and\/ $R_2$\/ are  Lagrangian submanifolds. By an
arbitrarily small perturbation of\/ $\omega_{f}$,  we can obtain
another symplectic form on $E(2)$ with respect to which\/ $F,S$\/
are still symplectic and $R_1$\/ and/or\/  $R_2$\/ are also
symplectic.
\end{lemma}

\begin{proof}
There is a symplectic form\/ $\omega_{f}$\/ on\/ $E(2)$\/ which is
induced by the elliptic fibration\/ $f\! :\! E(2)\rightarrow
\mathbb{CP}^1$, essentially as the sum of symplectic forms in the
fiber and the base (see \cite{thurston}). With respect to
$\omega_{f}$ a regular fiber $F$ and section $S$ are symplectic,
whereas the rim tori $R_1$ and $R_2$ are Lagrangian since the
circles $C_1$ and $C_2$ lie in $F$ and $\partial D^2$ is embedded
in a section. Since each $[R_i]$ is non-torsion and in fact
$[R_1]$ and $[R_2]$ are linearly independent, as a consequence of
the following more general lemma, we know that $\omega_{f}$ could
be slightly perturbed in order to make $R_1$ and/or $R_2$
symplectic.
\end{proof}

\begin{lemma}[cf.$\;$Lemma 1.6 in \cite{g:sum}]\label{lemma:lagrangian
perturbation}
Let $X$ be a closed\/ $4$-manifold with a symplectic form $\omega$
with respect to which
closed, connected and disjoint submanifolds
$\Sigma_1, \Sigma_2, \dots, \Sigma_r$
are Lagrangian. Suppose that the homology classes
$[\Sigma_1],[\Sigma_2], \dots, [\Sigma_r]$
are non-torsion and linearly independent. Then there exists an arbitrarily
small perturbation
$\omega'$ of $\omega$ which is symplectic and with respect to which all
surfaces
$\Sigma_1, \Sigma_2, \dots, \Sigma_r$ are symplectic submanifolds.
\end{lemma}

To prove the above lemma,
one needs to choose a closed 2-form $\Omega$ on $X$\/ such that
$\int_{\Sigma_i}
\Omega > 0\,$ for each $i$. Then\/
$\omega':= \omega'_s := \omega + s \Omega\,$
is a suitable perturbation for sufficiently small constant\/ $s > 0$.

\begin{theorem}\label{theorem:homology}
Fix a pair of integers\/ $q\geq 2$\/ and\/ $m\geq 1$.
\\
{\rm (i)}\/ The embedded torus\/ $S^1_x \times \phi_{q,m}(\hat{B})
\subset E(2)$ is a symplectic submanifold with respect to the
symplectic form\/ $\omega_f$, and represents the homology class\/
$q[F]+m[R_1]$.
\\
{\rm (ii)}\/ The embedded torus\/ $S^1_x \times
\psi_{q,m}(\hat{B}) \subset E(2)$\/ represents\/ $m[F]+q[R_1]$,
and there is a symplectic form on $E(2)$ with respect to which
this torus is a symplectic submanifold.
\end{theorem}

\begin{proof}
(i)\/\/ Without loss of generality, we may assume that the
restriction of\/ $\omega_f$ to the subset\/ $Y=\partial(\nu
F)\times[0,1]$\/ is given by (\ref{eq:symplectic form}).  This
immediately implies that\/ $S^1_x \times \phi_{q,m}(\hat{B})$\/
embeds symplectically into $E(2)$.  The link surgery gluing data
$\mathfrak{D}$\/ in (\ref{eq:data}) of Lemma~\ref{lemma:E(2)}
directly gives the homology class of\/ $S^1_x \times
\phi_{q,m}(\hat{B})$\/ in $E(2)$\/ since we have\/
$[\phi_{q,m}(\hat{B})] = q [\mu (K)] + m [\lambda (K)] \in H_1(S^3
\setminus \nu L; \zz)$, and $S^1_x$ gets identified with $C_1$.

(ii)\/\/ In Section~\ref{sec:cylinder}, we have already shown that
the torus $S^1_x \times \psi_{q,m}(\hat{B})$\/ is a symplectic
submanifold of\/ $Y = \partial(\nu F)\times[0,1] \hspace{1pt}$
with respect to the symplectic form $\omega_s$ for any\/ $s>0$.
By definition (\ref{def:omega}), $\omega_s = dx \wedge dy + r dr
\wedge d\theta + s \cdot dx \wedge d\theta$\/ near the boundary of
$Y$. Choosing the perturbation (which makes only $R_1$ symplectic)
in Lemma~\ref{lemma:symplectic form on E(2)} carefully (e.g.
$\Omega = dx \wedge d\theta$\/ with respect to the local
coordinates in which $\omega_{f}=dx \wedge dy + r dr \wedge
d\theta$) we could make sure that there exists a symplectic form
$\omega'$ on
$$ E(2) \;\cong\; [E(1) \setminus \nu F] \cup [\partial(\nu
F)\times[0,1]\hspace{1pt} ] \cup [E(1)\setminus \nu F] $$
which restricts
(up to isotopy)
to\/ $dx \wedge dy + r dr \wedge d\theta + s \cdot dx \wedge d\theta$\/ near
the boundary\/ $\partial
[\partial(\nu F)\times[0,1]\hspace{1pt} ]$. This allows us to extend
$\omega_s$ to the closed manifold
$E(2)$.

The link surgery gluing data $\mathfrak{D}$\/ in (\ref{eq:data})
again gives the homology class of\/ $S^1_x \times
\psi_{q,m}(\hat{B})$\/ in $E(2)$\/ since we have\/
$[\psi_{q,m}(\hat{B})] = m [\mu (K)] + q [\lambda (K)]$\/  this
time around.
\end{proof}

\begin{remark}\label{remark:switch}
Recall that we chose the factorization $F=C_1\times C_2$ in the
link surgery gluing data $\mathfrak{D}$\/ in (\ref{eq:data}) of
Lemma~\ref{lemma:E(2)}. If instead we had chosen the (reverse
order) identification\/ $F=C_2\times C_1$, then the tori\/ $S^1_x
\times \phi_{q,m}(\hat{B})$\/ and\/ $S^1_x \times
\psi_{q,m}(\hat{B})$\/ would have represented the homology
classes\/ $q[F] + m[R_2]\,$ and\/ $m[F] + q[R_2]\,$ in\/
$H_2(E(2);\zz)$, respectively.
\end{remark}


\section{Alexander Polynomials Corresponding to Particular Braids}
\label{sec:alexander}

\begin{figure}[!ht]
\begin{center}
\includegraphics[scale=.7]{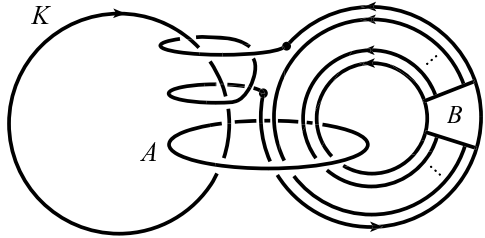}
\end{center}
\caption{3-component link $L \cup \psi_{q,m}(\hat{B})$\/ with
$m=2$} \label{fig:threelink}
\end{figure}

\begin{figure}[!ht]
\begin{center}
\includegraphics[scale=.6]{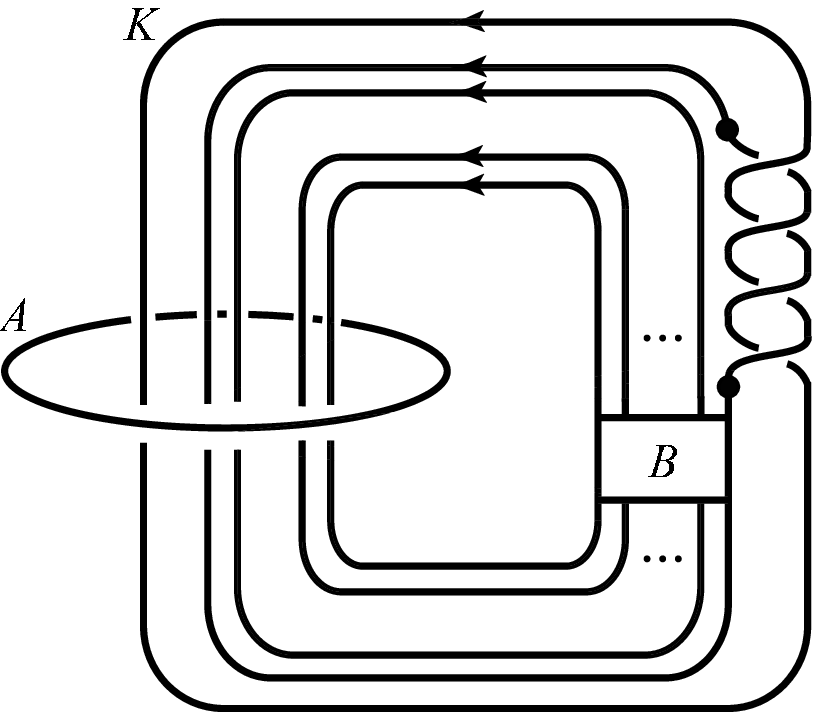}
\end{center}
\caption{3-component link $L \cup \psi_{q,m}(\hat{B})$\/ with
$m=2$\/ in braid form} \label{fig:braidform}
\end{figure}

In order to distinguish the isotopy classes of the homologous
symplectic tori we constructed in the previous section, we will
compute the Seiberg-Witten invariants of 4-manifolds that are
obtained as the fiber sum of $E(2)$ along these tori and the
rational elliptic surface $E(1)$ along one of its regular fibers.
We will see that the Seiberg-Witten invariant of such a 4-manifold
is essentially the Alexander polynomial of the 3-component link
obtained from the braid $B$\/ as seen in
Figures~\ref{fig:threelink} and \ref{fig:braidform}.  Both figures
are for the embedded tori\/ $S^1_x\times \psi_{q,m}(\hat{B})$, and
the corresponding pictures for\/ $S^1_x\times
\phi_{q,m}(\hat{B})$\/ are obtained by simply relabelling the
$K$\/ component $A$\/ and vice versa.

\begin{figure}[!ht]
\begin{center}
\includegraphics[scale=.6]{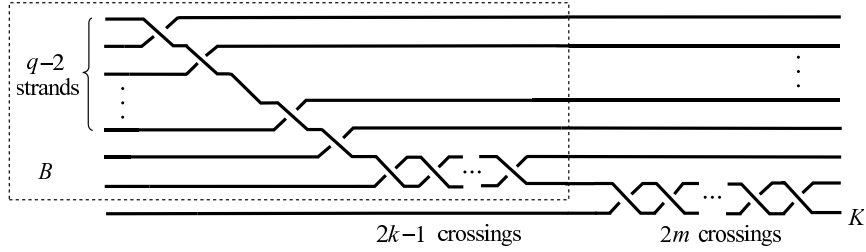}
\end{center}
\caption{$(q+1)$-strand braid $B(q;k,m)$\/ with
$k,m\geq 1$, $q\geq 2$}
\label{fig:simplebraid}
\end{figure}

In this section we will present the ``simplest'' family of braids
that is most amenable to the computation of the Alexander
polynomials of the corresponding links. A generic member
$B=B_{k,q}$ of this family is shown in
Figure~\ref{fig:simplebraid} as the upper left part (inside the
dotted rectangle) of the braid $B(q;k,m)$, for which the desired
3-component link $L \cup \psi_{q,m}(\hat{B})$ is $A \cup
\hat{B}(q;k,m)$, where $A$\/ is the axis of the closed braid
$\hat{B}(q;k,m)$ as well as one of the components of the Hopf link
$L=K \cup A$.

Similarly, we have $L\cup \phi_{q,m}(\hat{B})=K\cup
\hat{B}(q;k,m)$, where $K$\/ now denotes the axis of the braid
$B(q;k,m)$ and $A$\/ is the bottom strand in
Figure~\ref{fig:simplebraid}.

\begin{remark}\label{remark:q=0 or m=0}
Note that we are using the same family of braids $B$\/ as in
\cite{ep1}.  Consequently, when $m=0$ and $q\geq 2$, we obtain families of
tori representing either $q[F]$ or $q[R_i]$ that
we already constructed in \cite{ep1}.
\end{remark}

\begin{lemma}\label{lemma:alexander}
Let\/ $\Delta_{q;k,m} (x,s,t) $\/ denote the
Alexander polynomial of the three-component link\/ $L \cup
\psi_{q,m}(\hat{B}) = A \cup \hat{B}(q;k,m)\hspace{1pt}$, where
the variables\/ $x$, $s$ and\/ $t$\/ correspond to the axis\/ $A$,
unknot\/ $K$\/ and the closed braid\/ $\hat{B}$\/ respectively.
Then\/ $\Delta_{q;k,m} (x,s,t)= $
$$ 1-x(st)^m+ x \cdot
\frac{(xt)^{q-1}-1}{xt-1} \left[ t^{2k-1}+t(s-1)\,
\frac{t^{2k-1}+1}{t+1} \cdot \frac{(st)^m-1}{st-1}
-x(st)^mt^{2k-1}\right] . $$  The Alexander polynomial of the
link\/ $L \cup \phi_{q,m}(\hat{B}) = K \cup
\hat{B}(q;k,m)\hspace{1pt}$ is given by $\Delta_{q;k,m} (s,x,t)$,
i.e. the polynomial obtained from $\Delta_{q;k,m} (x,s,t)$ by
switching the variables $x$ and $s$.
\end{lemma}

\begin{proof}
The braid group on $q$\/ strands is generated by the elementary
braid transpositions\/ $\sigma_1,\dots ,
\sigma_{q-1}\hspace{1pt}$, where $\sigma_i$ denotes the crossing
of the (\hspace{1pt}$i$\hspace{1pt}+1)st strand over the $i$-th.
Note that

$$ B(q;k,m) \,=\: \sigma_{q}\,\sigma_{q-1} \,\cdots \,
\sigma_3 \,\sigma_2^{2k-1} \,\sigma_1^{2m}\, . $$

By Theorem~1 in \cite{morton}, we have
\begin{eqnarray}\label{eq:zero}
&& \Delta_{q;k,m} \, := \, \Delta_{q;k,m} (x,s,t) \, = \\ && \det\left(I-x\, C^{(q)}_{q}(t)
C^{(q)}_{q-1}(t) \cdots\, C^{(q)}_3(t) [C^{(q)}_2(t)]^{2k-1}
[C^{(q)}_1(s)C^{(q)}_1(t)]^m \right)\,
, \nonumber
\end{eqnarray}
where $C^{(q)}_i(a)$ denotes the following\/ $q\times q$\/ matrix
which differs from the identity matrix $I$\/ only in the three
places shown on the $i$-th row.
\[ C^{(q)}_i (a) \; :=\;
\left( \begin{array}{ccccccc} 1&&&&&&\\
&\ddots&&&&&\\&&1&&&&\\
&&a&-a&1&&\\
&&&&1&&\\
&&&&&\ddots& \\
&&&&&&1
\end{array}\right)
\!\!\begin{array}{r}
\\
\\
\\
\\
\\
\\
\\ .
\\
\end{array} \]
When\/ $i=1$\/ or\/ $i=q\hspace{1pt}$, the matrix is truncated
appropriately to give two non-zero entries in row $i$.

The main step of this proof is showing that $ D_{q;k,m}= xt D_{q-1;k,m}$
for all $q\geq 2$, where
$D_{q;k,m}:= \Delta_{q+1;k,m}-\Delta_{q;k,m}$ and
$\Delta_{1;k,m} := 1- x (st)^m$.
During this process we get
\begin{eqnarray}\label{eq:ahalf}
\Delta_{2;k,m}\!\!\! &=& \! 1-x^2(st)^mt^{2k-1} \\ &+& \!\!
x \left[t^{2k-1} -(st)^m +(s-1)t \left( \frac{t^{2k-1}+1}{t+1} \right)
\left( \frac{(st)^m-1}{st-1} \right) \right] \,  . \nonumber
\end{eqnarray}
This calculation leads to
$$D_{1;k,m}=x\left[t^{2k-1} + (s-1)t\left(\frac{t^{2k-1}+1}{t+1}\right)
\left( \frac{(st)^m-1}{st-1} \right)\right] - x^2(st)^mt^{2k-1}$$
and hence $D_{q;k,m}=$
$$(xt)^{q-1} \left\{ x\left[t^{2k-1} +
(s-1)t\left(\frac{t^{2k-1}+1}{t+1}\right)
\left( \frac{(st)^m-1}{st-1} \right) \right]- x^2(st)^mt^{2k-1} \right\} .$$
By putting the pieces together we finish the proof of the lemma.

By Equation~(\ref{eq:zero}) $\Delta_{q+1;k,m}=
\det (I - x \Gamma_{q+1;k,m})$, where
$$\Gamma_{q+1;k,m} =  C^{(q+1)}_{q+1}(t)
C^{(q+1)}_{q}(t) \cdots\, C^{(q+1)}_3(t) \!\!
\left[C^{(q+1)}_2(t)\right]^{2k-1}
\left[C^{(q+1)}_1(s)C^{(q+1)}_1(t)\right]^m .$$
Note that
\vspace{-5 pt}
$$ C^{(q+1)}_i (t)\;=\;
\left( \begin{array}{cccc} &&&0\\
&C^{(q)}_i(t)&&\vdots\\
&&&0 \\
0&\dots&0&1
\end{array}\right) $$
for $i \in \{ 1,2, \dots , q-1 \},\;$ so we must have $\Gamma_{q+1;k,m} = $
\begin{eqnarray*}
&=& C^{(q+1)}_{q+1}(t) C^{(q+1)}_{q}(t)
\left( \begin{array}{cccc} &&&0\\
&C^{(q)}_{q}(t)&&\vdots\\
&&&0 \\
0&\dots&0&1
\end{array}\right)^{\hspace{-8pt}
\begin{array}{l}  -1 \end{array}}
\hspace{-5pt}\left( \begin{array}{cccc} &&&0\\
&\Gamma_{q;k,m}&&\vdots\\
&&&0 \\
0&\dots&0&1
\end{array}\right) \\[5 pt]
&=&
\left( \begin{array}{cccc} 1&&&\\
&\ddots&&\\
&&1&1 \\
&&t&0
\end{array}\right)
\left( \begin{array}{cccc} &&&0\\
&\Gamma_{q;k,m}&&\vdots\\
&&&0 \\
0&\dots&0&1
\end{array}\right)\!\!\begin{array}{r}
\\
\\
\\ .
\end{array}
\end{eqnarray*}
Hence it follows that
\begin{equation}\label{eq:one} \Gamma_{q+1;k,m} \;=\;
\left( \begin{array}{cc} &0\\
\Gamma_{q;k,m}&\vdots\\
&0\\
&1 \\
t \, (\Gamma_{q;k,m})_{( q\hspace{1pt} ,\,\ast\hspace{1pt} )}&0
\end{array}\right)
\end{equation}
and
\begin{equation}\label{eq:oneandahalf}
\; I-x \Gamma_{q+1;k,m} \:=\:
\left( \begin{array}{cc} &0\\
I-x \Gamma_{q;k,m}&\vdots\\
&0\\
&-x \\
-xt \, (\Gamma_{q;k,m})_{(q\hspace{1pt} ,\,\ast\hspace{1pt})} &1
\end{array}\right)
\!\!\begin{array}{r}
\\
\\
\\
\\ ,
\end{array}
\end{equation}
where\/ $(\Gamma_{q;k,m})_{(q\hspace{1pt} ,\,\ast\hspace{1pt})}$\/
denotes the last row of $\Gamma_{q;k,m}\hspace{1pt}$.
When we calculate the determinant of the matrix $I - x \Gamma_{q+1;k,m}$
by expanding along its last column we get the following equality:
\begin{eqnarray}\label{eq:two}
&& \; \det (I-x\Gamma_{q+1;k,m} )\; = \\
&& \det ( I-x\Gamma_{q;k,m}) - (-x) t \left[\det (I-x\Gamma_{q;k,m}) -
\det (I-x\Gamma_{q-1;k,m})\right]\, . \nonumber
\end{eqnarray}

To prove the above equality for $q \geq 3$, observe that, in this case,
all but the last row of
the minor of the matrix $I-x\Gamma_{q+1;k,m}$ corresponding to the
entry $-x$ in the last column are the same as the rows of
$I-x\Gamma_{q;k,m}\hspace{1pt}$, and the last row of the minor is
$t$\/ times the last row of $I-x\Gamma_{q;k,m}$ except for the last
entry. In the minor, this entry is $0$, whereas in
$I-x\Gamma_{q;k,m}\hspace{1pt}$ this entry is $1$ (since
Equation~(\ref{eq:one}) shows that the last diagonal entry of
$\Gamma_{q;k,m}$ is $0$ as long as $q \geq 3$). This observation is why the
determinant of the minor corresponding to $-x$ is $t$ times the
difference between the determinant of $I-x\Gamma_{q;k,m}$ and the
determinant of the minor of $I-x\Gamma_{q;k,m}$ obtained by deleting
the last row and the last column (and this minor is nothing but
$I-x\Gamma_{q-1;k,m}$).

For $q=2$, Equation~(\ref{eq:two}) is proved
by direct calculation of $\Delta_{2;k,m}$ and $\Delta_{3;k,m}$.
Note that we {\it{defined}}\/ $\Delta_{1;k,m}$ to be $1-x(st)^m$. In fact,
once $I-x\Gamma_{2;k,m}$ is verified to be
$$\left(\!\! \begin{array}{cc}
 1-x (st)^m & x(s-1)\frac{(st)^m-1}{st-1} \\[7pt]
  -xt(st)^m\frac{t^{2k-1}+1}{t+1} &
\;\; 1+x(s-1)t\frac{t^{2k-1}+1}{t+1}\frac{(st)^m-1}{st-1}
+x t^{2k-1}
\end{array} \!\!\right)$$
with the help of the equalities
$$\left( \begin{array}{cc}
  1 & 0 \\
  t & -t
\end{array} \right)^{2k-1} = \;\left(\! \begin{array}{cc}
  1 & 0 \\ t\cdot\frac{t^{2k-1}+1}{t+1} & -t^{2k-1}
\end{array} \!\right) $$
and
$$ \left(\! \begin{array}{cc}
  st & -s+1 \\
  0 & 1
\end{array} \!\right)^m =
\;\left(\!\! \begin{array}{cc}
  (st)^m & -(s-1)\frac{(st)^m-1}{st-1} \\
  0 & 1
\end{array} \!\!\right) ,$$
one not only gets Equation~(\ref{eq:ahalf}) regarding $\Delta_{2;k,m}$,
but also the matrix $I-x\Gamma_{3;k,m}$ by using
Equation~(\ref{eq:oneandahalf}). As a result of expanding the $3 \times 3$
matrix
$I-x\Gamma_{3;k,m}$ along its last column, it is easily seen that
$$\Delta_{3;k,m}=\:\Delta_{2;k,m}+ x \{t[\Delta_{2;k,m}-(1-x(st)^m)]\} \, .$$

Equation~(\ref{eq:two}) and the calculations above give
\begin{eqnarray*}
 D_{q;k,m}\!\!\! &=& \!\!\! xt D_{q-1;k,m} \\
&=&\!\!\! (xt)^{q-1} x
\left[ t^{2k-1}+(s-1)t \,\frac{t^{2k-1}+1}{t+1} \cdot \frac{(st)^m-1}{st-1}
-x(st)^mt^{2k-1} \right] \, .
\end{eqnarray*}
Finally,
the formula in the statement of the lemma
is a consequence of
\[
\Delta_{q;k,m} \,=\: \Delta_{1;k,m}+D_{1;k,m}+ \cdots + D_{q-1;k,m} \]
\[ \hspace{1.25cm} =\: 1-x(st)^m+  \,\frac{(xt)^{q-1}-1}{xt-1} \,
D_{1;k,m}\: .
\qedhere
\]
\end{proof}

\begin{corollary}\label{corollary:terms}
The number of nonzero terms in the polynomial $\Delta_{q;k,m}(x,s,t)$
is equal to\/
$[6-4q]+k[2(m+1)(q-1)]$.
\end{corollary}
\begin{proof}
The polynomial $\Delta_{q;k,m}$ could be written as
\begin{eqnarray*}
\Delta_{q;k,m} &=& 1+ x\{-(st)^m+P_{k,m}(s,t)\} \\
&&+\: x^2\{t[-(st)^mt^{2k-2}+P_{k,m}(s,t)]\} \\
&&+\: \cdots \\
&&+\: x^{q-1}\{t^{q-2}[-(st)^mt^{2k-2}+P_{k,m}(s,t)]\} \\
&&+\: x^q\{-(st)^mt^{2k+q-3}\}\: ,
\end{eqnarray*}
where $$P_{k,m}(s,t)= t^{2k-1}+(s-1)t\,\frac{(st)^m-1}{st-1}\cdot
\frac{t^{2k-1}+1}{t+1}
\ .$$
A direct count of nonzero terms in $P_{k,m}(s,t)$ gives $2km+2k-3$,
and as a consequence, for $0<i<q$,
the number of nonzero terms in $\Delta_{q;k,m}$ that are divisible by $x^i$
but not divisible
by $x^{i+1}$ is $2km+2k-4$. The formula is then easily obtained
as a result of an effort to write the desired
expression in a form that emphasizes the dependence of the count on
$k$ when $m$ and $q$ are fixed.
\end{proof}


\section{Non-Isotopy:  Seiberg-Witten Invariants}
\label{sec:sw}

In Section~\ref{sec:tori in E(2)}, for each $i \in \{1,2\}$, $m
\geq 1$ and $q \geq 2$ we explained the construction of a
symplectic torus representing\/ $q[F]+m[R_i]$\/ or\/
$m[F]+q[R_i]$\/ using a suitable $q$-component braid $B$.  Let
$\jmath$\/ denote either $\phi_{q,m}$ or\/ $\psi_{q,m}$.  The
4-manifold $E(2) \#_{T=F} E(1)$, obtained as the fiber sum of
$E(1)$\/ along a regular fiber $F$\/ with $E(2)$ along one of
these tori\/ $T\! :=S^1_x \times \jmath(\hat{B})$\/ we
constructed, is easily seen to be diffeomorphic to the link
surgery manifold $(L \cup \jmath(\hat{B})) (\mathfrak{D}')$, where
$\mathfrak{D}'$ is the link surgery gluing data
$$ \left(\{ (\mu(K),\lambda(K)),(\lambda(A),-\mu(A)),
(\lambda(\jmath(\hat{B})),-\mu(\jmath(\hat{B})))\}, \{ E(1),
F=C_1\times C_2 \}_{i=1}^3 \right) .$$ In
Section~\ref{sec:alexander}, we looked at a particular family of
braids $B=B_{k,q}$ for which
$$ L_{\jmath}\, :=\, L \cup \jmath (\hat{B})
\;= \left\{ \!\! \begin{array}{lcl}  K \cup
\hat{B}(q;k,m) & {\rm  if } & \jmath=\phi_{q,m}\, ,\\[5pt]
A \cup \hat{B}(q;k,m) & {\rm  if } & \jmath=\psi_{q,m} \, .
\end{array}\right.
$$
In this section, we will distinguish the symplectic tori that
come from this family of braids by comparing the Seiberg-Witten
invariants of $L_{\jmath}(\mathfrak{D}')$.

Recall that the Seiberg-Witten invariant\/
$\overline{SW}_{\!\!X}$\/ of a 4-manifold $X$\/ (satisfying\/
$b_2^+(X)>1$) can be thought of as an element of the group ring of
$H_2(X;\zz)$, i.e.\/ $\overline{SW}_{\!\!X} \in \zz [ H_2( X ;
\zz) ]$. If we write\/ $\overline{SW}_{\!\!X} = \sum_g a_g g
\hspace{1pt}$, then we say that\/ $g\in H_2( X ; \zz)$\/ is a
Seiberg-Witten \emph{basic class}\/ of $X$\/ if\/ $a_g\neq 0$.
Since the Seiberg-Witten invariant of a 4-manifold is a
diffeomorphism invariant, so is the total number of Seiberg-Witten
basic classes.

Regarding the Seiberg-Witten invariants of
$L_{\jmath}(\mathfrak{D}')$, we have the following lemma which is
an easy consequence of the gluing formulas for the Seiberg-Witten
invariant in \cite{fs:knots}, \cite{doug:pft3} and
\cite{Taubes:T^3}. Detailed arguments can be found in \cite{ep1},
\cite{McMullen-Taubes} or \cite{vidussi:smooth}.

\begin{lemma}\label{lemma:seiberg-witten}
Let\/ $\iota:[S^1\times(S^3\setminus\nu L_{\jmath} )]\rightarrow
L_{\jmath}(\mathfrak{D}')$\/ be the inclusion map. Let\/
$\xi:=\iota_{\ast}[S^1\times\mu(A)],$
$\tau:=\iota_{\ast}[S^1\times\mu(K)],$
$\zeta:=\iota_{\ast}[S^1\times\mu(\jmath(\hat{B}))]\in
H_2(L_{\jmath}(\mathfrak{D}') ;\zz).$  Then the Seiberg-Witten
invariant of\/ $L_{\jmath}(\mathfrak{D}')$ is
\begin{eqnarray*}
\overline{SW}_{\! L_{\jmath}(\mathfrak{D}')} \!\!&=&\!\!
\Delta^{{\rm sym}}_{L_{\jmath}}(\xi^2,\tau^2,\zeta^2) \\[5pt]
&=&\!\! \left\{ \!\!\begin{array}{ccc}
\tau^{-q}\xi^{-m}\zeta^{-(2k+q+m-3)}\Delta_{q;k,m}(\tau^2,\xi^2,\zeta^2)
&{\rm if}& \jmath=\phi_{q,m} \, ,
\\[9pt]
\xi^{-q}\tau^{-m}\zeta^{-(2k+q+m-3)}\Delta_{q;k,m}(\xi^2,\tau^2,\zeta^2)
&{\rm if}& \jmath=\psi_{q,m} \, ,
\end{array}\right.
\end{eqnarray*}
where\/ $\Delta_{q;k,m}$ is the Alexander polynomial in Lemma\/
$\ref{lemma:alexander}$, and\/ $\Delta^{{\rm sym}}$ stands for the
symmetrized Alexander polynomial.
\end{lemma}

Note that\/ $\xi,\tau$ and $\zeta$ are linearly independent in
$H_2(L_{\jmath}(\mathfrak{D}');\zz)$ as in Proposition 3.2 of
\cite{McMullen-Taubes}. As a consequence of
Corollary~\ref{corollary:terms}, the number of Seiberg-Witten
basic classes of $L_{\jmath}(\mathfrak{D}')$ depends on $k$\/ for
fixed $q\ge 2$ and $m\ge 1$. Hence, for fixed triple $q$, $m$ and
$\jmath\hspace{1pt}$, the family of 4-manifolds\/
$\{L_{\jmath}(\mathfrak{D}')\}_{k\geq 1}$ are all pairwise
non-diffeomorphic. On the other hand, the diffeomorphism type of
$L_{\jmath}(\mathfrak{D}') \cong E(2) \#_{T=F} E(1)$ only depends
on the isotopy type of\/ $T$.

This finishes the proof of Theorem~\ref{theorem:main}. In fact,
one can easily see that the tori we constructed are different even
under self-diffeomorphisms of $E(2)$.


\section{Generalization to Other Symplectic 4-Manifolds}
\label{sec:generalization}

For certain elliptic surfaces, our result easily generalizes.
Since our tori will remain
non-isotopic even after fiber sum and link surgery (cf.$\;$\cite{fs:ipam}),
we immediately obtain the analogue of Theorem~\ref{theorem:main} for
the fiber sums
$E(n)=E(2) \#_F E(n-2)$ for $n\geq 3$, and the knot surgery manifolds
$$E(n)_K := \: K(\{(\alpha_1, \beta_1)=(\mu(K),\lambda(K))\},
\{E(n),F=C_1\times C_2\})$$
for any fibred knot $K\subset S^3$ and $n\geq 2$.
(Note that the knot
$K$\/ needs to be fibred
to ensure that $E(n)_K$ is symplectic,
and
$E(n)_K$ can also be viewed as the fiber sum\/
$E(n-1)\#_F E(1)_K\hspace{1pt}$.)
Also note that
an infinite subset
of our homologous symplectic tori will continue to remain
different under self-diffeomorphisms of these
symplectic 4-manifolds, since the number of Seiberg-Witten basic classes
of the corresponding link surgery
manifolds always goes to infinity as $k\rightarrow \infty$ and $q,m$
are fixed.

In particular, we recover and
generalize Vidussi's result (Corollary 1.2 in \cite{vidussi:lagrangian})
on the non-isotopic
symplectic representatives of primitive homology classes on certain
knot surgery
manifolds $E(2)_K$ (also see \cite{fs:lagrangian}).

For more general symplectic 4-manifolds, note that the Hopf link
will give us any fiber sum manifold like $E(2)$. More precisely,
if $Z$\/ is obtained as the symplectic fiber sum along symplectic
tori of self-intersection $0$, then by choosing a suitable link
surgery gluing data, we can symplectically embed $S^1_x \times
\jmath(\hat{B})$ in $Z$. In order to distinguish these tori we can
still use Seiberg-Witten theory, but we need some extra
assumptions to make use of the gluing formulas for the
Seiberg-Witten invariant.

\begin{theorem}\label{thm:generalization}
Suppose that $F_i$ is a symplectically embedded\/ $2$-torus in a
closed symplectic\/ $4$-manifold $Z_i$ with $b^+_2(Z_i) > 1$,
$[F_i]^2=0$ and $H^1(Z_i \setminus \nu F_i; \zz)=0$, for each $i
\in \{1,2 \}$. Let $Z=Z_1\#_{F_1=F_2}Z_2$\/ be the symplectic
fiber sum of $Z_1$ and $Z_2$ along $F_1$ and $F_2$.  Let\/ $[F]$
and $[R]$ be the homology classes of $F_1=F_2$ and a rim torus in
$Z$, respectively. Then for any pair of positive integers
$(q,m)\neq (1,1)$ there exists an
infinite family of pairwise non-isotopic symplectic tori
representing the homology class\/ $q[F]+m[R]\in H_2(Z;\zz)$.
\end{theorem}

\begin{proof}
Let $L$\/ denote a Hopf link in $S^3$ as before.
We can express $Z=L(\mathfrak{D}'')$, where
\[   \mathfrak{D}'' :=\: \big( \{ (\mu(K),\lambda(K)),(\lambda(A),-\mu(A))\},
\{ Z_i, F_i=C^i_1\times C^i_2 \}_{i=1}^2 \big)\, . \] The rim
torus $R$\/ in the lemma is given by the Cartesian product $C^i_1
\times
\partial D^2$, where $D^2$ is a normal disk in $\nu F_i\cong
F_i\times D^2$. We need to compute the Seiberg-Witten invariants
of the corresponding link surgery manifolds
$L_{\jmath}(\mathfrak{D}''')$, where
\begin{eqnarray*}
\mathfrak{D}'''  \!\!\! &:=& \!\!\! \big(\{ (\mu(K),\lambda(K)),
(\lambda(A),-\mu(A)),
(\lambda(\jmath(\hat{B})),-\mu(\jmath(\hat{B}))) \}, \\
&& \{ Z_i, F_i=C^i_1\times C^i_2 \}_{i=1}^2  \cup \{ E(1),
F=C_1\times C_2 \} \big)\, .
\end{eqnarray*}
Just as in \cite{ep1}, the assumption that
$H^1(Z_i \setminus \nu F_i; \zz)=0$ ($i=1,2$) is crucial.
It allows us to conclude that the homology classes\/ $[F]$ and $[R]$\/
are linearly independent in $H_2(Z;\zz)$ as in Proposition 3.2
of \cite{McMullen-Taubes}.
It also implies that the relative Seiberg-Witten invariants are
\[
\overline{SW}_{\! Z_i\setminus \nu F_i} =\:
([F_i]^{-1}-[F_i])\cdot\overline{SW}_{\! Z_i} \:\neq 0 \,
\]
by Corollary 20 in \cite{doug:pft3}.  Hence the Seiberg-Witten
invariants of $L_{\jmath}(\mathfrak{D}''')$\/ can be computed
using the standard gluing formulas as before. The rest of the
proof is the same as the proof of Theorem~\ref{theorem:main}. Once
again, to conclude that there are infinitely many tori that remain
different under self-diffeomorphisms of $Z$, we observe that, for
fixed pair $q$\/ and $m$, the number of Seiberg-Witten basic
classes of $L_{\jmath}(\mathfrak{D}''')$\/ goes to infinity as
$k\rightarrow \infty$. Non-isotopy is more simply obtained from a
homology basis argument due to Fintushel and Stern
(cf.$\;$\cite{fs:ipam}).
\end{proof}

\begin{remark}
The conclusion of Theorem~\ref{thm:generalization} may still apply
even when $b_2^+(Z_i)=1$. In that case, one must take care and
define $\overline{SW}_{\!\!Z_i}:=\overline{SW}_{\!\!Z_i,
F_i}^{\hspace{1pt}\pm}$\/ (see \cite{fs:knots} and
\cite{doug:pft3}). In general, for a closed 4-manifold $X$ with
$b_2^+(X)=1$, it is not automatic that\/ $\overline{SW}_{\!\!X}$
is a finite sum and $\overline{SW}_{\!\!X}\neq 0$\/ for a
symplectic $X$. If indeed $\overline{SW}_{\!\!Z_i}\neq 0$\/ and is
a finite sum, then Theorem~\ref{thm:generalization} will still be
valid for such $Z$. However if\/ $\overline{SW}_{\!\!Z_i}=0$\/ or
is an infinite sum, then there seems to be no systematic method
currently available to check whether the tori in our family are
mutually non-isotopic in $Z$ or not. An ad hoc method for a
particularly simple infinite sum case is presented in
\cite{ep:E(1)_K} for a slightly different family of tori
(corresponding to embeddings $\phi_{1,m}$).
\end{remark}


\smallskip
\subsection*{Acknowledgments}
We would like to thank Ronald Fintushel, Ian Hambleton,
Maung Min-Oo, Sa\v{s}o Strle
and Stefano Vidussi for their encouragement and helpful comments.
The figures were produced by the second author using
Adobe$^{\circledR}$
{\sl Illustrator}$\hspace{1pt}^{\circledR}$
Version 10. Some computations in
Section~\ref{sec:alexander} were verified with the aid of {\sl
Maple}$\hspace{1pt}^{\circledR}$   Version 8.

\smallskip



\begin{thebibliography}{[ADK]}

\bibitem[ADK]{adk} D. Auroux, S.K. Donaldson and L. Katzarkov:
Luttinger surgery along Lagrangian tori and non-isotopy for
singular symplectic plane curves, {\it Math. Ann.} {\bf 326}
(2003), 185--203.

\bibitem[EP1]{ep1}  T. Etg\"u and B.D. Park:
Non-isotopic symplectic tori in the same homology class, {\sl
preprint}.  Available at arXiv:math.GT/0212356.

\bibitem[EP2]{ep:E(1)_K}  T. Etg\"u and B.D. Park:
Homologous non-isotopic
symplectic tori in homotopy rational elliptic surfaces, {\sl
preprint}.  Available at arXiv:math.GT/0307029.

\bibitem[EP3]{ep2}  T. Etg\"u and B.D. Park:
Homologous non-isotopic Lagrangian tori in symplectic 4-manifolds,
{\sl in preparation}.

\bibitem[FS1]{fs:knots} R. Fintushel and R.J. Stern: Knots,
links and $4$-manifolds, {\it Invent. Math.} {\bf 134} (1998),
363--400.

\bibitem[FS2]{fs:non-isotopic} R. Fintushel and R.J. Stern:
Symplectic surfaces in a fixed homology class, {\it J. Differential Geom.}
{\bf 52} (1999), 203--222.

\bibitem[FS3]{fs:lagrangian} R. Fintushel and R.J. Stern:
Invariants for Lagrangian tori, {\sl preprint}.
Available at arXiv:math.SG/0304402.

\bibitem[FS4]{fs:ipam}  R. Fintushel and R.J. Stern:
Tori in symplectic $4$-manifolds, {\sl preprint}.


\bibitem[Go]{g:sum} R.E. Gompf: A new construction of symplectic manifolds,
{\it Ann. of Math.} {\bf 142} (1995), 527--595.

\bibitem[GS]{gs}  R.E. Gompf and A.I. Stipsicz:
{\sl $4$-Manifolds and Kirby Calculus\/}, Graduate Studies in
Mathematics {\bf 20}, Amer. Math. Soc., 1999.

\bibitem[MT]{McMullen-Taubes} C.T. McMullen and C.H. Taubes:
$4$-manifolds with inequivalent symplectic forms and $3$-manifolds
with inequivalent fibrations, {\it Math. Res. Lett.\/} {\bf 6}
(1999), 681--696.

\bibitem[Mo]{morton}  H.R. Morton:
The multivariable Alexander polynomial for a closed braid, in {\sl
Low-dimensional Topology}, ed. Hanna Nencka, {\it Contemporary
Mathematics} {\bf 233}, Amer. Math. Soc. (1999), 167--172. Also
available at arXiv:math.GT/9803138.

\bibitem[Pa]{doug:pft3}  B.D. Park:  A gluing formula for
the Seiberg-Witten invariant along $T^3$, {\it Michigan Math.
J.\/} {\bf 50} (2002), 593--611.

\bibitem[Ta]{Taubes:T^3}  C.H. Taubes:
The Seiberg-Witten invariants and $4$-manifolds with essential
tori, {\it Geom. Topol.\/} {\bf 5} (2001),  441--519.

\bibitem[Th]{thurston} W.P. Thurston:
Some simple examples of symplectic manifolds,
{\it Proc. Amer. Math. Soc.\/} {\bf 55} (1976), 467--468.

\bibitem[V1]{vidussi:smooth}  S. Vidussi:
Smooth structure of some symplectic surfaces, {\it Michigan Math.
J.} {\bf 49} (2001), 325--330.

\bibitem[V2]{vidussi:non-isotopic} S. Vidussi:
Nonisotopic symplectic tori in the fiber class of elliptic surfaces,
{\sl preprint}.  Available at {\tt
http://www.math.ksu.edu/\~{}vidussi/}

\bibitem[V3]{vidussi:lagrangian}  S. Vidussi:
Lagrangian surfaces in a fixed homology class:
Existence of knotted Lagrangian tori,  {\sl preprint}.
Available at {\tt
http://www.math.ksu.edu/\~{}vidussi/}

\bibitem[V4]{vidussi:E(1)_K}  S. Vidussi:
Symplectic tori in homotopy $E(1)$'s,  {\sl preprint}.  Available
at\\ {\tt http://www.math.ksu.edu/\~{}vidussi/}


\end{thebibliography}
\end{document}